\documentclass[11pt,a4paper]{amsart}
\usepackage[utf8]{inputenc}
\usepackage{amsmath,amssymb, amsthm, amsfonts, gensymb, commath, float, mathtools, enumerate, breqn, mathdots}
\usepackage{tabularx}
\usepackage{multirow}
\usepackage{esvect}
\usepackage{subcaption}
\usepackage{tikz}
\usetikzlibrary{positioning, arrows, arrows.meta, cd}
\usetikzlibrary{shapes,backgrounds,shapes,positioning,petri,topaths}

\theoremstyle{theorem}
\newtheorem{theorem}{Theorem}[section]

\newtheorem{conjecture}[theorem]{Conjecture}

\def\beq#1#2\eeq{%
        \begin{equation}%
        \label{#1}%
            #2%
        \end{equation}%
   }

\usepackage{color}
\usepackage[dvipsnames]{xcolor}

\definecolor{red}{RGB}{216, 27, 96}
\definecolor{blue}{RGB}{30, 136, 229}
\definecolor{green}{RGB}{0, 77, 64}
\definecolor{yellow}{RGB}{255, 193, 7}
\usepackage{graphicx}
\usepackage{epstopdf}
\usepackage[text={15cm,23cm}]{geometry} 
\usepackage[colorlinks=true,%
            linkcolor=red!50!black,%
            citecolor=blue!50!black,%
            urlcolor=darkgray]{hyperref}  

\DeclareMathOperator{\Conv}{Conv}

\theoremstyle{definition}
\newtheorem{example}[theorem]{Example}
\newtheorem{definition}[theorem]{Definition}

\title[Markov Polynomials]{Saturation of Markov Polynomials}

\author{S.J. Evans}
\address{Department of Mathematical Sciences,
Loughborough University, Loughborough LE11 3TU, UK}
\email{S.J.Evans@lboro.ac.uk}

\begin{document}

\begin{abstract}
    Solutions to the Markov equation appear in many mathematical contexts. We aim to build on the understanding of them by proving a recent conjecture about Markov polynomials; solutions to a generalised version of the Markov equation. The proof we provide is a constructive argument based on the Markov snake graph, a combinatorial object related to Markov numbers, deepening the connection between the Markov equation and combinatorics.
\end{abstract}

\maketitle

\section{Introduction}
Markov numbers (solutions to the Markov equation) play a major role in various areas of mathematics \cite{Markov}, \cite{Aig}. Markov polynomials are the Laurent polynomial solutions to the generalised Markov equation
\begin{equation}
\label{eqn:Gen Markov}
    X^2 + Y^2 + Z^2 = k(x, y, z) XYZ, \qquad k(x, y, z) = \frac{x^2 + y^2 + z^2}{xyz},
\end{equation}
in the three variables $x, y, z$. They have previously been studied in the context of cluster algebras, which is where it was proven that these polynomials are Laurent (see Fomin and Zelevinsky \cite{FZ} and also Propp \cite{Pro}).

We studied the terms that appear in these polynomials more closely in \cite{EVW2}. 
Firstly, the Markov numbers can be indexed by the rational numbers in $[0, 1]$ (from Frobenius \cite{Frob}). Given a rational $\rho \in [0, 1]$, we define $m_\rho$ to be the corresponding Markov number. We can calculate the Markov number for a given rational in several ways, including using combinatorics which we will see in more detail. In \cite{EVW2} we prove that given a rational number $\rho = \frac{a}{b}$, the corresponding Markov polynomial, $M_\rho$ has the form
\begin{equation*}
    M_{a/b}(x, y, z) = \frac{P_{a/b}(x^2, y^2, z^2)}{x^{a-1} y^{b-1} z^{a+b-1}},
\end{equation*}
where $P_{a/b}(u, v, w)$ is homogeneous of total degree $a+b-1$ (scaling the even powers). The numerator can be written as
\begin{equation*}
    P_\rho(u, v, w) = \sum A_{ij} u^i v^j w^{a+b-1-i-j}.
\end{equation*}
The Newton polygon, $\Delta_\rho$ is defined to be
\begin{equation}
    \Delta_\rho = \Delta(P_\rho) := \Conv \lbrace (i, j) : A_{ij} \neq 0 \rbrace
\end{equation}
and in \cite{EVW2} we prove the following theorem.

\begin{theorem}[Theorem 3.2 in \cite{EVW2}]
    \label{thm:Newton Poly}
    Given a rational $\rho=\frac{a}{b}$, the Newton polygon $\Delta_\rho$ is the area on the $ij$-plane with $i, j \geq 0$ satisfying the conditions
    \begin{equation}
    \label{eqn:Newton Poly}
         \left\{
        \begin{aligned}
            & \frac{i}{a} + \frac{j}{b} \geq 1 \\
            & i + j \leq a+b-1.
        \end{aligned}
        \right.
    \end{equation}
\end{theorem}

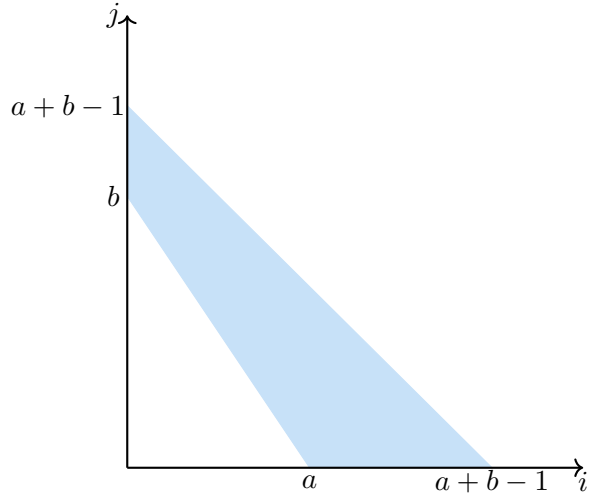
\begin{figure}[H]
\begin{center}
\begin{tikzpicture}  [scale=0.6]
        \filldraw[blue!25] (8, 0) -- (0, 8) -- (0, 6) -- (4, 0) -- (8, 0);
        
        \draw[thick, ->] (0, 0) -- (0, 10);
        \draw[thick, ->] (0, 0) -- (10, 0);
        
        \node at (8, -0.3) {$a+b-1$};
        \node at (4, -0.3) {$a$};
        \node at (-1.3, 8) {$a+b-1$};
        \node at (-0.3, 6) {$b$};

        \node at (-0.3, 10) {$j$};
        \node at (10, -0.3) {$i$};
\end{tikzpicture}
\end{center}
\end{figure}

We further claim that all of the integer lattice points in the set $\Delta_\rho \cap \mathbb{Z}^2$ do in fact appear as terms in the Markov polynomial (i.e., the corresponding monomials all have non-zero coefficients). This is the Saturation conjecture \cite{EVW2}.

\begin{conjecture}[Conjecture 4.3 in \cite{EVW2}]
\label{conj:Saturation}
    The terms that appear in the numerator of a Markov polynomial $M_{\rho}$ are precisely those corresponding to the integer lattice points in the Newton polygon $\Delta_{\rho}$. 
\end{conjecture}

In this paper we use the known association between Markov polynomials and the combinatorial object of snake graphs \cite{Pro} in order to provide a constructive proof of the Saturation conjecture.

\section{Combinatorial Interpretation of Markov}
It is known that Markov numbers, $m_\rho$ can be determined by counting the number of perfect matchings of the so-called Markov snake graph. This snake graph can be constructed from the rational $\rho$ and appears in various works \cite{Aig}, \cite{Cohn}, \cite{Pro}, \cite{Reu}.

Given a rational number $\rho = \frac{a}{b} \in [0, 1]$, draw the line with gradient $\rho$ from the origin to the point $(b, a)$. We then look at the closest path between integer lattice points below this diagonal (known as the Christoffel path). Here, we show an example, where $\rho = \frac{3}{5}$ with the Christoffel path given in blue.
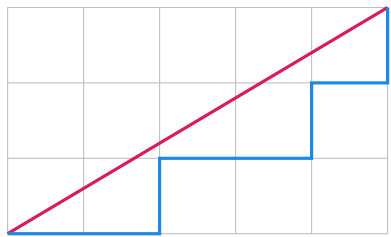
\begin{figure}[H]
\begin{center}
\begin{tikzpicture}
    \draw[lightgray] (0, 0) grid (5, 3);
    \draw[red, very thick] (0, 0) -- (5, 3);
    \draw[blue, very thick] (0, 0) -- (2, 0) -- (2, 1) -- (4, 1) -- (4, 2) -- (5, 2) -- (5, 3);
\end{tikzpicture}
\end{center}
\end{figure}
The Markov snake graph is then formed by placing half-unit boxes above the the Christoffel path, omitting the first and final box.
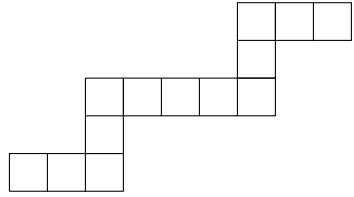
\begin{figure}[H]
\begin{center}
\begin{tikzpicture}[scale=0.5]
    \draw (1, 0) grid (4, 1);
    \draw (3, 1) grid (4, 3);
    \draw (4, 2) grid (8, 3);
    \draw (7, 3) grid (8, 5);
    \draw (8, 4) grid (10, 5);
\end{tikzpicture}
\end{center}
\end{figure}
The number of perfect matchings of this graph is then the Markov number $m_\rho$.

This idea was extended to the generalised case by adding weights to these edges to form the weighted Markov snake graph \cite{IMPV}. Each perfect matching had a monomial associated to it (as the product of the edge-weights involved in the particular perfect matching). The (numerator of the) corresponding Markov polynomial was proven to be the sum of these monomials \cite{Pro}.

This simplest explanation for the rules of adding weights is to consider, moving from left-to-right, every other box. The vertical edges of these boxes have weight $x$ and the horizontal edges have weight $y$. The remaining edges of the whole snake graph that remain unlabelled have weight $z$. Continuing the example above:

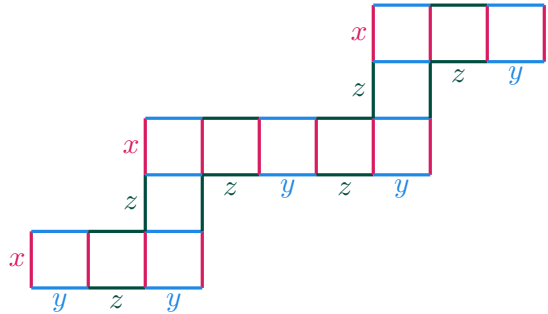
\begin{figure}[H]
\begin{center}
\begin{tikzpicture}[scale=0.75]
    \draw[red, very thick] (0, 0) -- (0, 1);
    \draw[red, very thick] (1, 0) -- (1, 1);
    \draw[blue, very thick] (0, 0) -- (1, 0);
    \draw[blue, very thick] (0, 1) -- (1, 1);
    \draw[green, very thick] (1, 0) -- (2, 0);
    \draw[green, very thick] (1, 1) -- (2, 1);

    \draw[red, very thick] (2, 0) -- (2, 1);
    \draw[red, very thick] (3, 0) -- (3, 1);
    \draw[blue, very thick] (2, 0) -- (3, 0);
    \draw[blue, very thick] (2, 1) -- (3, 1);
    \draw[green, very thick] (2, 1) -- (2, 2);
    \draw[green, very thick] (3, 1) -- (3, 2);

    \draw[red, very thick] (2, 2) -- (2, 3);
    \draw[red, very thick] (3, 2) -- (3, 3);
    \draw[blue, very thick] (2, 2) -- (3, 2);
    \draw[blue, very thick] (2, 3) -- (3, 3);
    \draw[green, very thick] (3, 2) -- (4, 2);
    \draw[green, very thick] (3, 3) -- (4, 3);

    \draw[red, very thick] (4, 2) -- (4, 3);
    \draw[red, very thick] (5, 2) -- (5, 3);
    \draw[blue, very thick] (4, 2) -- (5, 2);
    \draw[blue, very thick] (4, 3) -- (5, 3);
    \draw[green, very thick] (5, 2) -- (6, 2);
    \draw[green, very thick] (5, 3) -- (6, 3);
    
    \draw[red, very thick] (6, 2) -- (6, 3);
    \draw[red, very thick] (7, 2) -- (7, 3);
    \draw[blue, very thick] (6, 2) -- (7, 2);
    \draw[blue, very thick] (6, 3) -- (7, 3);
    \draw[green, very thick] (6, 3) -- (6, 4);
    \draw[green, very thick] (7, 3) -- (7, 4);

    \draw[red, very thick] (6, 4) -- (6, 5);
    \draw[red, very thick] (7, 4) -- (7, 5);
    \draw[blue, very thick] (6, 4) -- (7, 4);
    \draw[blue, very thick] (6, 5) -- (7, 5);
    \draw[green, very thick] (7, 4) -- (8, 4);
    \draw[green, very thick] (7, 5) -- (8, 5);

    \draw[red, very thick] (8, 4) -- (8, 5);
    \draw[red, very thick] (9, 4) -- (9, 5);
    \draw[blue, very thick] (8, 4) -- (9, 4);
    \draw[blue, very thick] (8, 5) -- (9, 5);

    \node[red, very thick] at (-0.25, 0.5) {$x$};
    \node[green, very thick] at (1.75, 1.5) {$z$};
    \node[red, very thick] at (1.75, 2.5) {$x$};
    \node[green, very thick] at (5.75, 3.5) {$z$};
    \node[red, very thick] at (5.75, 4.5) {$x$};
    \node[blue, very thick] at (0.5, -0.25) {$y$};
    \node[green, very thick] at (1.5, -0.25) {$z$};
    \node[blue, very thick] at (2.5, -0.25) {$y$};
    \node[green, very thick] at (3.5, 1.75) {$z$};
    \node[blue, very thick] at (4.5, 1.75) {$y$};
    \node[green, very thick] at (5.5, 1.75) {$z$};
    \node[blue, very thick] at (6.5, 1.75) {$y$};
    \node[green, very thick] at (7.5, 3.75) {$z$};
    \node[blue, very thick] at (8.5, 3.75) {$y$};
\end{tikzpicture}
\end{center}
\end{figure}

For clarity, only the outermost edges are labelled. However, all edges are weighted, with those of the same colour being the same weight ($x$ for red edges, $y$ for blue and $z$ for green).

\section{Proof Strategy}
To prove saturation, it is enough to show that there exists at least one perfect matching corresponding to each lattice point on the Newton polygon.

\subsection{Newton Polygon Anatomy}
First we define some terms regarding the anatomy of the Newton polygon.
\begin{definition}
We refer to the lattice points on the Newton polygon based on the diagonal line $i + j = c$ (constant) that they lie on.
\begin{itemize}
    \item {\it Full diagonals} are the diagonals for which the point $(0, c)$ lies on the boundary of the Newton polygon, i.e., $b \leq c \leq a+b-1$.
    \item {\it Partial diagonals} are the diagonals for which the point $(0, c)$ (and potentially other points in the first quadrant on the diagonal) lies outside of the Newton polygon, i.e., $c < b$.
    \item {\it Left-most point} is the point on a diagonal with minimal $i$-coordinate that lies within the Newton polygon. Note that for full diagonals this will be $(0, c)$ and for partial diagonals it will be $(i, c-i)$ for some $i \geq 1$.
\end{itemize}
\end{definition}

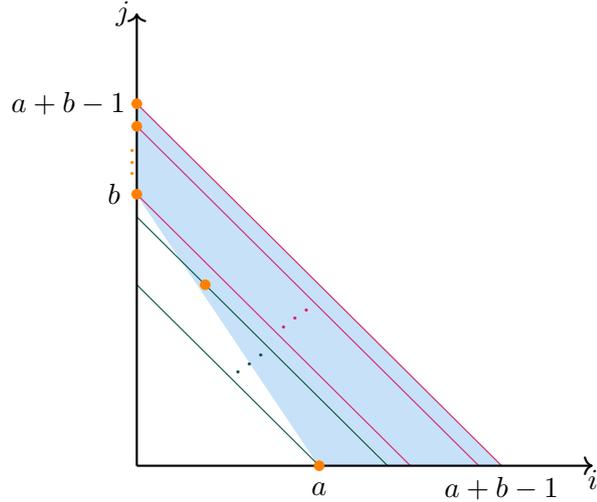
\begin{figure}[H]
\begin{center}
\begin{tikzpicture}[scale=0.6]
        \filldraw[blue!25] (8, 0) -- (0, 8) -- (0, 6) -- (4, 0) -- (8, 0);
        \draw[red] (8, 0) -- (0, 8);
        \draw[red] (7.5, 0) -- (0, 7.5);
        \node[red] at (3.45, 3.45) {$\iddots$};
        \draw[red] (6, 0) -- (0, 6);
        \draw[green] (5.5, 0) -- (0, 5.5);
        \node[green] at (2.45, 2.45) {$\iddots$};
        \draw[green] (4, 0) -- (0, 4);
        
        \draw[thick, ->] (0, 0) -- (0, 10);
        \draw[thick, ->] (0, 0) -- (10, 0);

        \filldraw[orange] (0, 8) circle (3pt);
        \filldraw[orange] (0, 7.5) circle (3pt);
        \node[orange] at (-0.1, 6.9) {$\vdots$};
        \filldraw[orange] (0, 6) circle (3pt);
        \filldraw[orange] (1.5, 4) circle (3pt);
        \filldraw[orange] (4, 0) circle (3pt);
        
        \node at (8, -0.5) {$a+b-1$};
        \node at (4, -0.5) {$a$};
        \node at (-1.5, 8) {$a+b-1$};
        \node at (-0.5, 6) {$b$};

        \node at (-0.3, 10) {$j$};
        \node at (10, -0.3) {$i$};
\end{tikzpicture}
\caption{Newton polygon illustrating the full diagonals (red), partial diagonals (green) and left-most points (orange).}
\label{fig: Newton polygon anatomy}
\end{center}
\end{figure}

\subsection{Path}

We first claim that we can reduce this problem to just finding a perfect matching corresponding to the left-most point on each diagonal $i+j = c$ inside the Newton polygon. Then to move along each diagonal we use a set of {\it Traversing Operations}.

So we have to find perfect matchings corresponding to the left-most points within the Newton polygon on each of the diagonals
\begin{equation*}
    i + j = a + b - 1,\ a + b - 2,\ \dots,\ b,\ \dots,\ a.
\end{equation*}
We produce these by starting at an {\it Initial Perfect Matching (IPM)} and performing local operations on this matching (from left-to-right) to move down to left-most point of each subsequent diagonal.

The way in which we construct these perfect matchings also ensures that we do not `miss' points, and that we do indeed hit the first relevant point on each diagonal.

\section{Constructive Algorithm}
\subsection{Initial Perfect Matching}
We define the {\it Initial Perfect Matching (IPM)} as the perfect matching consisting of all $y$-weight edges. That is, we take the horizontal edges of every other box. In the example above
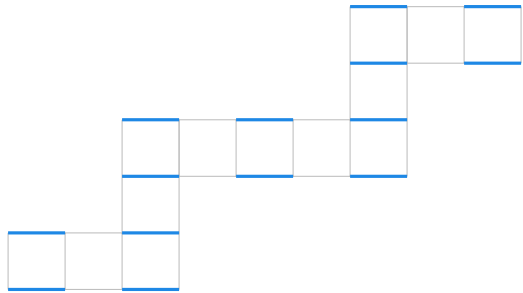
\begin{figure}[H]
\begin{center}
\begin{tikzpicture}[scale=0.75]
    \draw[lightgray] (0, 0) grid (3, 1);
    \draw[lightgray] (2, 1) grid (3, 3);
    \draw[lightgray] (3, 2) grid (7, 3);
    \draw[lightgray] (6, 3) grid (7, 5);
    \draw[lightgray] (7, 4) grid (9, 5);
    
    \draw[blue, very thick] (0, 0) -- (1, 0);
    \draw[blue, very thick] (0, 1) -- (1, 1);

    \draw[blue, very thick] (2, 0) -- (3, 0);
    \draw[blue, very thick] (2, 1) -- (3, 1);
 
    \draw[blue, very thick] (2, 2) -- (3, 2);
    \draw[blue, very thick] (2, 3) -- (3, 3);

    \draw[blue, very thick] (4, 2) -- (5, 2);
    \draw[blue, very thick] (4, 3) -- (5, 3);

    \draw[blue, very thick] (6, 2) -- (7, 2);
    \draw[blue, very thick] (6, 3) -- (7, 3);

    \draw[blue, very thick] (6, 4) -- (7, 4);
    \draw[blue, very thick] (6, 5) -- (7, 5);

    \draw[blue, very thick] (8, 4) -- (9, 4);
    \draw[blue, very thick] (8, 5) -- (9, 5);

\end{tikzpicture}
\end{center}
\end{figure}
Note that the number of these boxes will always be $a + b - 1$, due to the way in which we construct the snake graph.

Therefore, this perfect matching corresponds to the monomial $y^{2(a+b-1)}$ or the lattice point $(0, a+b-1)$, which is the first point on the first diagonal.

\subsection{Full Diagonals}
Now we construct the perfect matchings corresponding to the left-most point on each of the remaining full diagonals, i.e.,
\begin{equation*}
    i + j = a + b - 2,\ \dots,\ b+1,\ b.
\end{equation*}
First note that there are precisely $a-1$ of these diagonals, and $a-1$ regions of the IPM of the form
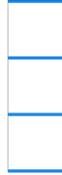
\begin{figure}[H]
\begin{center}
\begin{tikzpicture}[scale=0.75]
    \draw[lightgray] (0, 0) grid (1, 3);
    
    \draw[blue, very thick] (0, 0) -- (1, 0);
    \draw[blue, very thick] (0, 1) -- (1, 1);
    \draw[blue, very thick] (0, 2) -- (1, 2);
    \draw[blue, very thick] (0, 3) -- (1, 3);
\end{tikzpicture}
\end{center}
\end{figure}
To move down $k \leq a-1$ diagonals we replace $k$ of these regions by the alternative perfect matching
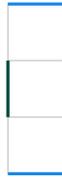
\begin{figure}[H]
\begin{center}
\begin{tikzpicture}[scale=0.75]
    \draw[lightgray] (0, 0) grid (1, 3);
    
    \draw[blue, very thick] (0, 0) -- (1, 0);
    \draw[green, very thick] (0, 1) -- (0, 2);
    \draw[green, very thick] (1, 1) -- (1, 2);
    \draw[blue, very thick] (0, 3) -- (1, 3);
\end{tikzpicture}
\end{center}
\end{figure}
Each time we perform this operation on the perfect matching, we are replacing one copy of $y^2$ with $z^2$. Thus, for each one we move one point down in the lattice. Hence the first points on each of these full diagonals is
\begin{equation*}
    (0, a+b-2),\ (0, a+b-3),\ \dots,\ (0, b+1),\ (0, b),
\end{equation*}
as expected.

\subsection{Partial Diagonals}
For the perfect matching on the diagonal $i + j = b$, all of the vertical regions of the snake graph of the IPM have been adapted to incorporate $z$-edges. Continuing the example above, we have
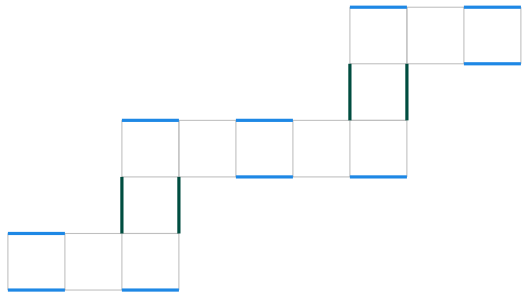
\begin{figure}[H]
\begin{center}
\begin{tikzpicture}[scale=0.75]
    \draw[lightgray] (0, 0) grid (3, 1);
    \draw[lightgray] (2, 1) grid (3, 3);
    \draw[lightgray] (3, 2) grid (7, 3);
    \draw[lightgray] (6, 3) grid (7, 5);
    \draw[lightgray] (7, 4) grid (9, 5);
    
    \draw[blue, very thick] (0, 0) -- (1, 0);
    \draw[blue, very thick] (0, 1) -- (1, 1);

    \draw[blue, very thick] (2, 0) -- (3, 0);
    \draw[green, very thick] (2, 1) -- (2, 2);
    \draw[green, very thick] (3, 1) -- (3, 2);
    \draw[blue, very thick] (2, 3) -- (3, 3);

    \draw[blue, very thick] (4, 2) -- (5, 2);
    \draw[blue, very thick] (4, 3) -- (5, 3);

    \draw[blue, very thick] (6, 2) -- (7, 2);
    \draw[green, very thick] (6, 3) -- (6, 4);
    \draw[green, very thick] (7, 3) -- (7, 4);
    \draw[blue, very thick] (6, 5) -- (7, 5);

    \draw[blue, very thick] (8, 4) -- (9, 4);
    \draw[blue, very thick] (8, 5) -- (9, 5);

\end{tikzpicture}
\end{center}
\end{figure}

We now perform a series of local operations on the perfect matching that move us from left-to-right. The operations are determined precisely by the shape of the snake graph (there are no choices to be made).

\subsubsection{First Operation}
To begin with, we are in one of two cases

\noindent {\bf Case 1 - Addition, $\oplus$:}
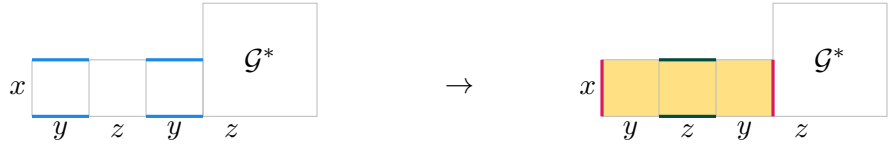
\begin{figure}[H]
\begin{center}
\begin{tikzpicture}[scale = 0.75]
    \draw[lightgray] (0, 0) grid (3, 1);
    \draw[lightgray] (3, 0) rectangle (5, 2);
    \draw[blue, very thick] (0, 0) -- (1, 0);
    \draw[blue, very thick] (0, 1) -- (1, 1);
    \draw[blue, very thick] (2, 0) -- (3, 0);
    \draw[blue, very thick] (2, 1) -- (3, 1);

    \node at (-0.25, 0.5) {$x$};
    \node at (0.5, -0.25) {$y$};
    \node at (1.5, -0.25) {$z$};
    \node at (2.5, -0.25) {$y$};
    \node at (3.5, -0.25) {$z$};
    \node at (4, 1) {$\mathcal{G}^*$};

    \begin{scope}[xshift=5cm]
        \node at (2.5, 0.5) {$\rightarrow$};
    \end{scope}
    
    \begin{scope}[xshift=10cm]
        \filldraw[yellow!50] (0, 0) rectangle (3, 1);
        \draw[lightgray] (0, 0) grid (3, 1);
        \draw[lightgray] (3, 0) rectangle (5, 2);
        \draw[red, very thick] (0, 0) -- (0, 1);
        \draw[green, very thick] (1, 0) -- (2, 0);
        \draw[green, very thick] (1, 1) -- (2, 1);
        \draw[red, very thick] (3, 0) -- (3, 1);

        \node at (-0.25, 0.5) {$x$};
        \node at (0.5, -0.25) {$y$};
        \node at (1.5, -0.25) {$z$};
        \node at (2.5, -0.25) {$y$};
        \node at (3.5, -0.25) {$z$};
        \node at (4, 1) {$\mathcal{G}^*$};
    \end{scope}
\end{tikzpicture}
\end{center}
\end{figure}
This represents a shift $(+1, -2)$ on the integer lattice.

\noindent {\bf Case 2 - Initial Step, $\Gamma_i^*$:}
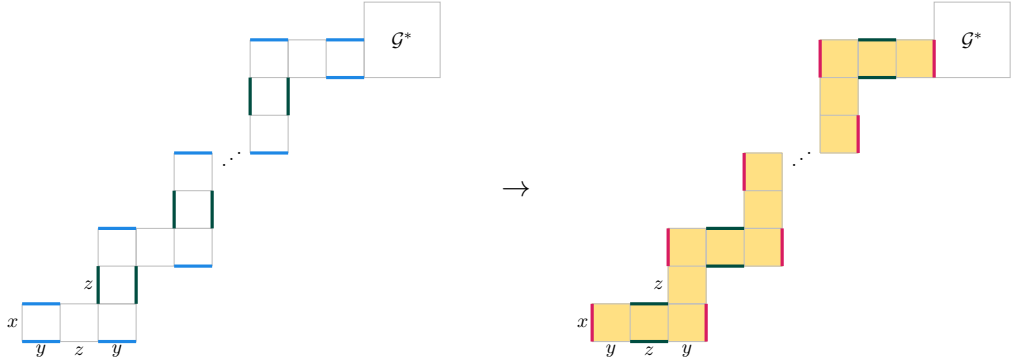
\begin{figure}[H]
\begin{center}
\begin{tikzpicture}[scale = 0.5, every node/.style={scale = 0.67}]
    \draw[lightgray] (-1, 0) grid (2, 1);
    \draw[lightgray] (1, 1) grid (2, 3);
    \draw[lightgray] (2, 2) grid (4, 3);
    \draw[lightgray] (3, 3) grid (4, 5);
    \draw[lightgray] (5, 5) grid (6, 8);
    \draw[lightgray] (6, 7) grid (8, 8);
    \draw[lightgray] (8, 7) rectangle (10, 9);
    \draw[blue, very thick] (-1, 0) -- (0, 0);
    \draw[blue, very thick] (-1, 1) -- (0, 1);
    \draw[blue, very thick] (1, 0) -- (2, 0);
    \draw[green, very thick] (1, 1) -- (1, 2);
    \draw[green, very thick] (2, 1) -- (2, 2);
    \draw[blue, very thick] (1, 3) -- (2, 3);
    \draw[blue, very thick] (3, 2) -- (4, 2);
    \draw[green, very thick] (3, 3) -- (3, 4);
    \draw[green, very thick] (4, 3) -- (4, 4);
    \draw[blue, very thick] (3, 5) -- (4, 5);
    
    \draw[blue, very thick] (5, 5) -- (6, 5);
    \draw[green, very thick] (5, 6) -- (5, 7);
    \draw[green, very thick] (6, 6) -- (6, 7);
    \draw[blue, very thick] (5, 8) -- (6, 8);
    \draw[blue, very thick] (7, 7) -- (8, 7);
    \draw[blue, very thick] (7, 8) -- (8, 8);

    \node at (-1.25, 0.5) {$x$};
    \node at (-0.5, -0.25) {$y$};
    \node at (0.5, -0.25) {$z$};
    \node at (1.5, -0.25) {$y$};
    \node at (0.75, 1.5) {$z$};
    \node at (4.5, 5) {$\iddots$};
    \node at (9, 8) {$\mathcal{G}^*$};

    \begin{scope}[xshift = 7.5cm, every node/.style={scale = 1}]
        \node at (4.5, 4) {$\rightarrow$};
    \end{scope}

    \begin{scope}[xshift = 15cm]
        \filldraw[yellow!50] (-1, 0) -- (2, 0) -- (2, 2) -- (4, 2) -- (4, 5) -- (3, 5) -- (3, 3) -- (1, 3) -- (1, 1) -- (-1, 1) -- (-1, 0);
        \filldraw[yellow!50] (5, 5) -- (6, 5) -- (6, 7) -- (8, 7) -- (8, 8) -- (5, 8) -- (5, 5);
        \draw[lightgray] (-1, 0) grid (2, 1);
        \draw[lightgray] (1, 1) grid (2, 3);
        \draw[lightgray] (2, 2) grid (4, 3);
        \draw[lightgray] (3, 3) grid (4, 5);
        \draw[lightgray] (5, 5) grid (6, 8);
        \draw[lightgray] (6, 7) grid (8, 8);
        \draw[lightgray] (8, 7) rectangle (10, 9);
        \draw[red, very thick] (-1, 0) -- (-1, 1);
        \draw[green, very thick] (0, 0) -- (1, 0);
        \draw[green, very thick] (0, 1) -- (1, 1);
        \draw[red, very thick] (2, 0) -- (2, 1);
        \draw[red, very thick] (1, 2) -- (1, 3);
        \draw[green, very thick] (2, 2) -- (3, 2);
        \draw[green, very thick] (2, 3) -- (3, 3);
        \draw[red, very thick] (4, 2) -- (4, 3);
        \draw[red, very thick] (3, 4) -- (3, 5);

        \draw[red, very thick] (6, 5) -- (6, 6);
        \draw[red, very thick] (5, 7) -- (5, 8);
        \draw[green, very thick] (6, 7) -- (7, 7);
        \draw[green, very thick] (6, 8) -- (7, 8);
        \draw[red, very thick] (8, 7) -- (8, 8);
    
        \node at (-1.25, 0.5) {$x$};
        \node at (-0.5, -0.25) {$y$};
        \node at (0.5, -0.25) {$z$};
        \node at (1.5, -0.25) {$y$};
        \node at (0.75, 1.5) {$z$};
        \node at (4.5, 5) {$\iddots$};
        \node at (9, 8) {$\mathcal{G}^*$};
    \end{scope}
\end{tikzpicture}
\end{center}
\end{figure}
Here, $i$ denotes the number of `steps' that we move up. This operation represents a shift $(+(i+1), -(i+2))$ on the integer lattice. Note that $i = 0$ is effectively the addition operation.

\subsubsection{Remaining Operations}
At this point, we are again in one of two cases

\noindent {\bf Case 1 - Extension, $\otimes$:}
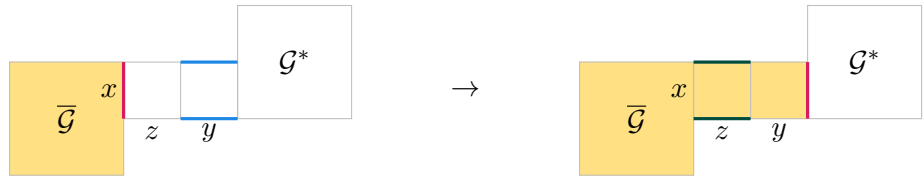
\begin{figure}[H]
\begin{center}
\begin{tikzpicture}[scale = 0.75]
    \filldraw[yellow!50] (-2, -1) rectangle (0, 1);
    \draw[lightgray] (0, 0) grid (2, 1);
    \draw[lightgray] (-2, -1) rectangle (0, 1);
    \draw[lightgray] (2, 0) rectangle (4, 2);
    \draw[red, very thick] (0, 0) -- (0, 1);
    \draw[blue, very thick] (1, 0) -- (2, 0);
    \draw[blue, very thick] (1, 1) -- (2, 1);

    \node at (-0.25, 0.5) {$x$};
    \node at (0.5, -0.25) {$z$};
    \node at (1.5, -0.25) {$y$};
    \node at (-1, 0) {$\overline{\mathcal{G}}$};
    \node at (3, 1) {$\mathcal{G}^*$};

    \begin{scope}[xshift=5cm]
         \node at (1, 0.5) {$\rightarrow$};
    \end{scope}

    \begin{scope}[xshift=10cm]
        \filldraw[yellow!50] (-2, -1) rectangle (0, 1);
        \filldraw[yellow!50] (0, 0) rectangle (2, 1);
        \draw[lightgray] (0, 0) grid (2, 1);
        \draw[lightgray] (-2, -1) rectangle (0, 1);
        \draw[lightgray] (2, 0) rectangle (4, 2);
        \draw[green, very thick] (0, 0) -- (1, 0);
        \draw[green, very thick] (0, 1) -- (1, 1);
        \draw[red, very thick] (2, 0) -- (2, 1);
    
        \node at (-0.25, 0.5) {$x$};
        \node at (0.5, -0.25) {$z$};
        \node at (1.5, -0.25) {$y$};
        \node at (-1, 0) {$\overline{\mathcal{G}}$};
        \node at (3, 1) {$\mathcal{G}^*$};
    \end{scope}
\end{tikzpicture}
\end{center}
\end{figure}
This represents a shift $(0, -1)$ on the integer lattice.

\noindent{\bf Case 2 - Step, $\Gamma_i$:}
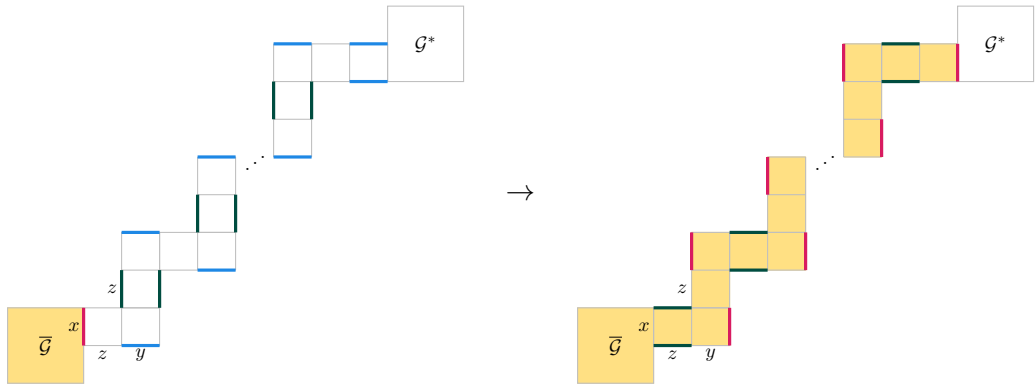
\begin{figure}[H]
\begin{center}
\begin{tikzpicture}[scale = 0.5, every node/.style={scale = 0.67}]
    \filldraw[yellow!50] (-2, -1) rectangle (0, 1);
    \draw[lightgray] (-2, -1) rectangle (0, 1);
    \draw[lightgray] (0, 0) grid (2, 1);
    \draw[lightgray] (1, 1) grid (2, 3);
    \draw[lightgray] (2, 2) grid (4, 3);
    \draw[lightgray] (3, 3) grid (4, 5);
    \draw[lightgray] (5, 5) grid (6, 8);
    \draw[lightgray] (6, 7) grid (8, 8);
    \draw[lightgray] (8, 7) rectangle (10, 9);
    \draw[red, very thick] (0, 0) -- (0, 1);
    \draw[blue, very thick] (1, 0) -- (2, 0);
    \draw[green, very thick] (1, 1) -- (1, 2);
    \draw[green, very thick] (2, 1) -- (2, 2);
    \draw[blue, very thick] (1, 3) -- (2, 3);
    \draw[blue, very thick] (3, 2) -- (4, 2);
    \draw[green, very thick] (3, 3) -- (3, 4);
    \draw[green, very thick] (4, 3) -- (4, 4);
    \draw[blue, very thick] (3, 5) -- (4, 5);
    
    \draw[blue, very thick] (5, 5) -- (6, 5);
    \draw[green, very thick] (5, 6) -- (5, 7);
    \draw[green, very thick] (6, 6) -- (6, 7);
    \draw[blue, very thick] (5, 8) -- (6, 8);
    \draw[blue, very thick] (7, 7) -- (8, 7);
    \draw[blue, very thick] (7, 8) -- (8, 8);

    \node at (-0.25, 0.5) {$x$};
    \node at (0.5, -0.25) {$z$};
    \node at (1.5, -0.25) {$y$};
    \node at (0.75, 1.5) {$z$};
    \node at (4.5, 5) {$\iddots$};
    \node at (9, 8) {$\mathcal{G}^*$};
    \node at (-1, 0) {$\overline{\mathcal{G}}$};

    \begin{scope}[xshift = 7.5cm, every node/.style={scale = 1}]
        \node at (4, 4) {$\rightarrow$};
    \end{scope}

    \begin{scope}[xshift = 15cm]
        \filldraw[yellow!50] (01, 0) -- (2, 0) -- (2, 2) -- (4, 2) -- (4, 5) -- (3, 5) -- (3, 3) -- (1, 3) -- (1, 1) -- (0, 1) -- (0, 0);
        \filldraw[yellow!50] (5, 5) -- (6, 5) -- (6, 7) -- (8, 7) -- (8, 8) -- (5, 8) -- (5, 5);
        \filldraw[yellow!50] (-2, -1) rectangle (0, 1);
        \draw[lightgray] (-2, -1) rectangle (0, 1);
        \draw[lightgray] (0, 0) grid (2, 1);
        \draw[lightgray] (1, 1) grid (2, 3);
        \draw[lightgray] (2, 2) grid (4, 3);
        \draw[lightgray] (3, 3) grid (4, 5);
        \draw[lightgray] (5, 5) grid (6, 8);
        \draw[lightgray] (6, 7) grid (8, 8);
        \draw[lightgray] (8, 7) rectangle (10, 9);
        \draw[green, very thick] (0, 0) -- (1, 0);
        \draw[green, very thick] (0, 1) -- (1, 1);
        \draw[red, very thick] (2, 0) -- (2, 1);
        \draw[red, very thick] (1, 2) -- (1, 3);
        \draw[green, very thick] (2, 2) -- (3, 2);
        \draw[green, very thick] (2, 3) -- (3, 3);
        \draw[red, very thick] (4, 2) -- (4, 3);
        \draw[red, very thick] (3, 4) -- (3, 5);

        \draw[red, very thick] (6, 5) -- (6, 6);
        \draw[red, very thick] (5, 7) -- (5, 8);
        \draw[green, very thick] (6, 7) -- (7, 7);
        \draw[green, very thick] (6, 8) -- (7, 8);
        \draw[red, very thick] (8, 7) -- (8, 8);
    
        \node at (-0.25, 0.5) {$x$};
        \node at (0.5, -0.25) {$z$};
        \node at (1.5, -0.25) {$y$};
        \node at (0.75, 1.5) {$z$};
        \node at (4.5, 5) {$\iddots$};
        \node at (9, 8) {$\mathcal{G}^*$};
        \node at (-1, 0) {$\overline{\mathcal{G}}$};
    \end{scope}
\end{tikzpicture}
\end{center}
\end{figure}
Again, $i$ denotes the number of `steps' that we move up. This operation represents a shift $(+i, -(i+1))$ on the integer lattice. Note that $i = 0$ is effectively the extension operation.

We perform these operations until $\mathcal{G}^*$ is empty.

\subsection{Traversing Diagonals}
Suppose we are trying to obtain a perfect matching corresponding to an arbitrary lattice point on the Newton polygon. To do this, we first locate the left-most point on the appropriate diagonal by terminating the process described above once the appropriate diagonal has been reached. Then we continue performing local operations from left-to-right along the snake graph but with a different pair of {\it traversal operations}, see below. We are again in one of the two cases from before but now the operations are different.

\noindent {\bf Case 1 - Twist, $\#$:}
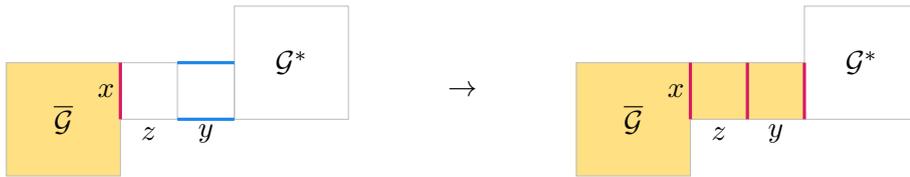
\begin{figure}[H]
\begin{center}
\begin{tikzpicture}[scale = 0.75]
    \filldraw[yellow!50] (-2, -1) rectangle (0, 1);
    \draw[lightgray] (0, 0) grid (2, 1);
    \draw[lightgray] (-2, -1) rectangle (0, 1);
    \draw[lightgray] (2, 0) rectangle (4, 2);
    \draw[red, very thick] (0, 0) -- (0, 1);
    \draw[blue, very thick] (1, 0) -- (2, 0);
    \draw[blue, very thick] (1, 1) -- (2, 1);

    \node at (-0.25, 0.5) {$x$};
    \node at (0.5, -0.25) {$z$};
    \node at (1.5, -0.25) {$y$};
    \node at (-1, 0) {$\overline{\mathcal{G}}$};
    \node at (3, 1) {$\mathcal{G}^*$};

    \begin{scope}[xshift=5cm]
         \node at (1, 0.5) {$\rightarrow$};
    \end{scope}

    \begin{scope}[xshift=10cm]
        \filldraw[yellow!50] (-2, -1) rectangle (0, 1);
        \filldraw[yellow!50] (0, 0) rectangle (2, 1);
        \draw[lightgray] (0, 0) grid (2, 1);
        \draw[lightgray] (-2, -1) rectangle (0, 1);
        \draw[lightgray] (2, 0) rectangle (4, 2);
        \draw[red, very thick] (0, 0) -- (0, 1);
        \draw[red, very thick] (1, 0) -- (1, 1);
        \draw[red, very thick] (2, 0) -- (2, 1);
    
        \node at (-0.25, 0.5) {$x$};
        \node at (0.5, -0.25) {$z$};
        \node at (1.5, -0.25) {$y$};
        \node at (-1, 0) {$\overline{\mathcal{G}}$};
        \node at (3, 1) {$\mathcal{G}^*$};
    \end{scope}
\end{tikzpicture}
\end{center}
\end{figure}

\noindent {\bf Case 2 - Pull-back, $\lhd$:}
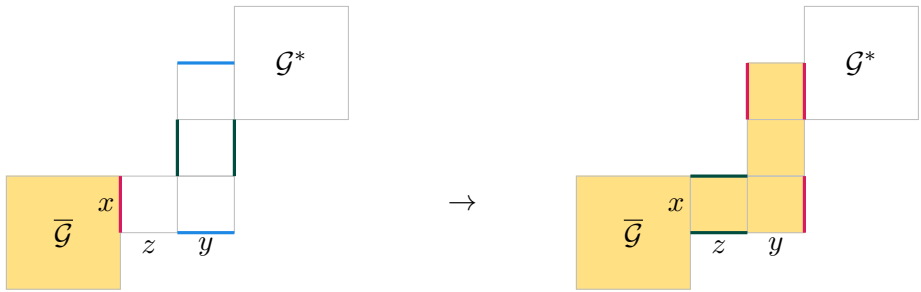
\begin{figure}[H]
\begin{center}
\begin{tikzpicture}[scale = 0.75]
    \filldraw[yellow!50] (-2, -1) rectangle (0, 1);
    \draw[lightgray] (0, 0) grid (2, 1);
    \draw[lightgray] (1, 1) grid (2, 3);
    \draw[lightgray] (-2, -1) rectangle (0, 1);
    \draw[lightgray] (2, 2) rectangle (4, 4);
    \draw[red, very thick] (0, 0) -- (0, 1);
    \draw[blue, very thick] (1, 0) -- (2, 0);
    \draw[green, very thick] (1, 1) -- (1, 2);
    \draw[green, very thick] (2, 1) -- (2, 2);
    \draw[blue, very thick] (1, 3) -- (2, 3);

    \node at (-0.25, 0.5) {$x$};
    \node at (0.5, -0.25) {$z$};
    \node at (1.5, -0.25) {$y$};
    \node at (-1, 0) {$\overline{\mathcal{G}}$};
    \node at (3, 3) {$\mathcal{G}^*$};

    \begin{scope}[xshift=5cm]
         \node at (1, 0.5) {$\rightarrow$};
    \end{scope}

    \begin{scope}[xshift=10cm]
        \filldraw[yellow!50] (-2, -1) rectangle (0, 1);
        \filldraw[yellow!50] (0, 0) -- (2, 0) -- (2, 3) -- (1, 3) -- (1, 1) -- (0, 1) -- (0, 0);
        \draw[lightgray] (0, 0) grid (2, 1);
        \draw[lightgray] (1, 1) grid (2, 3);
        \draw[lightgray] (-2, -1) rectangle (0, 1);
        \draw[lightgray] (2, 2) rectangle (4, 4);
        \draw[green, very thick] (0, 0) -- (1, 0);
        \draw[green, very thick] (0, 1) -- (1, 1);
        \draw[red, very thick] (2, 0) -- (2, 1);
        \draw[red, very thick] (1, 2) -- (1, 3);
        \draw[red, very thick] (2, 2) -- (2, 3);
    
        \node at (-0.25, 0.5) {$x$};
        \node at (0.5, -0.25) {$z$};
        \node at (1.5, -0.25) {$y$};
        \node at (-1, 0) {$\overline{\mathcal{G}}$};
        \node at (3, 3) {$\mathcal{G}^*$};
    \end{scope}
\end{tikzpicture}
\end{center}
\end{figure}
Both of these operations represent a shift $(+1, -1)$ on the integer lattice (i.e., moving along the diagonal).

Further note that if $\overline{\mathcal{G}}$ is empty then we start by swapping the first box. Originally, we  will have two $y$-weight edges and they are replaced by the $x$-weight ones.

\section{Covering all Lattice Points}
\subsection{Left-most Points}
Whilst our construction guarantees that we move from one diagonal to the one below it after each operation, it may not be immediately clear that this is the left-most point on this diagonal. In theory, there could be a lattice point inside the Newton polygon that we have missed. In order to convince ourselves that this will not be the case we look at the Christoffel word, of which we refer to \cite{Evans} for a more detailed summary.

The Christoffel word, $\omega_\rho(\alpha, \beta)$ is the word produced by tracking the direction of the Christoffel path. That is, every time we move horizontally we put an $\alpha$ and every time we move vertically, a $\beta$. Returning to the example of $\rho = \frac{3}{5}$ we have
\begin{equation*}
    \omega_{\frac{3}{5}}(\alpha, \beta) = \alpha \alpha \beta \alpha \alpha \beta \alpha \beta.
\end{equation*}
Taking the change of variables $\alpha \beta \rightarrow B$ and then $\alpha \rightarrow A$ we obtain the modified word
\begin{equation*}
    \widetilde{\omega_{\frac{3}{5}}}(A, B) = ABABB.
\end{equation*}

First note that, a string of $B$'s indicates a step up (except for the very final $B$) in the snake graph and will therefore result in a step operation when we reach that point in the snake graph. Similarly, consecutive $A$'s correspond to the extension operation (or addition if at the very start). Therefore, we claim that the shifts on the Newton polygon can be described entirely in terms of the lengths of strings of $B$'s between each pair of $A$'s.

In our example, $\rho = \frac{3}{5}$, we have one $B$ between the first pair of $A$'s and so our first operation is a step $\Gamma_1^*$. Then we have two consecutive $B$'s after the final $A$ and so the second operation is a step $\Gamma_1$ (only one step as this last $B$ does not result in a step on the snake graph; due to the omission of the final box when constructing the Markov snake along the Christoffel path).

Finally, we claim that the left-most points on each diagonal can also be determined from this modified word and that in fact they match the lattice points obtained from these shifts as described above.

This is because of the crucial observation that the {\it critical triangle} of the Newton polygon, which is the region bounded by the lines \cite{EVW2}
\begin{equation*}
    i < a, \qquad  j < b, \qquad \frac{i}{a} + \frac{j}{b} > 1,
\end{equation*}
displays similarity with the diagonal used to construct the snake graph (a rotation of $90\degree$). 

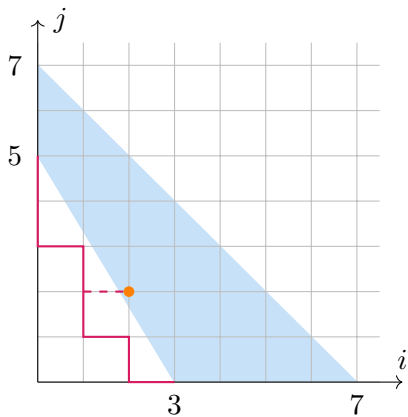
\begin{figure}[H]
\begin{center}
\begin{tikzpicture}
\begin{scope}[xshift = 8cm, scale = 0.6]
    \filldraw[blue!25] (7, 0) -- (0, 7) -- (0, 5) -- (3, 0) -- (7, 0);
    \draw[lightgray] (0, 0) grid (7.5, 7.5);
    
    \draw[->] (0, 0) -- (0, 8);
    \draw[->] (0, 0) -- (8, 0);
    \draw[thick, red] (0, 5) -- (0, 3) -- (1, 3) -- (1, 1) -- (2, 1) -- (2, 0) -- (3, 0);
    \draw[thick, red, dashed] (1, 2) -- (2, 2);
    \filldraw[orange] (2, 2) circle (3pt);
    
    \node at (7, -0.5) {$7$};
    \node at (3, -0.5) {$3$};
    \node at (-0.5, 7) {$7$};
    \node at (-0.5, 5) {$5$};

    \node at (0.5, 8) {$j$};
    \node at (8, 0.5) {$i$};
\end{scope}
\end{tikzpicture}
\caption{Diagram showing the Christoffel path superimposed on the Newton polygon.}
\label{fig:Superimposition}
\end{center}
\end{figure}

Therefore, the left-most points on the diagonal can be seen as the points at which the snake graph deviates from alternating between horizontal and vertical moves. At this point, if we continued to alternate between these two moves, we would cross the lower boundary of the Newton polygon and equivalently reach a lattice point that is on the appropriate diagonal and inside the Newton polygon. However, alternating horizontal and vertical moves corresponds precisely to $B$'s in the modified word. Therefore, $A$'s in the modified word indicate these left-most (first) points.

So, again, we are looking at the lengths of strings of $B$'s between two $A$'s. Let $i_j \geq 0$ denote the length of the string of $B$'s between the $(j-1)$th and $j$th $A$ in the word. Then, the first lattice point within the Newton polygon for each diagonal will be
\begin{equation*}
    (i_1 + 1, b - 2 - i_1), \quad (i_1 + i_2 + 1, b - 3 - i_1 - i_2), \quad \dots, \quad (a, 0).
\end{equation*}
The general form is
\begin{equation*}
    \left(1 + \sum_j i_j,\ b - (j+1) - \sum_j i_j\right).
\end{equation*}
It is clear that the shifts between the points on each diagonal is precisely the shift determined by the operation that we perform at this stage. This confirms that we do indeed obtain the left-most point on each diagonal.

\subsection{Traversing}
For traversing diagonals, we will always be able to reach the $i$-axis because, after all operations, the yellow shaded region no longer contains any $y$-weight edges. Hence, if we perform local operations all the way to the end of the graph then there will be no $y$ powers in the monomial corresponding to the whole perfect matching.

Note that for traversing the correspondence between the modified word and the operations is even more straightforward. A single $A$ implies a twist and a single $B$; a pull-back.

\section{Full Worked Example}
The process can be more clearly seen by example.
\begin{example}
    {\it Q. Produce the perfect matching corresponding to the point $(4, 2)$ on the Newton polygon $\Delta_\rho$ for $\rho = \frac{4}{7}$.}

    {\it A.} First we must locate the initial point on the diagonal 
    \begin{equation*}
        i + j = 4 + 2 = 6,
    \end{equation*}
    finding a corresponding perfect matching.
    
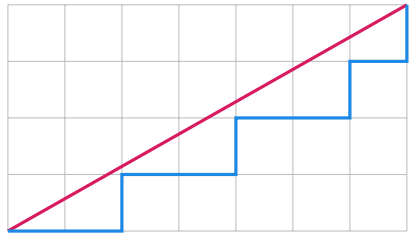
\begin{figure}[H]
    \begin{center}
    \begin{tikzpicture}
    \begin{scope}[scale = 0.75]
        \draw[lightgray] (0, 0) grid (7, 4);
        \draw[red, very thick] (0, 0) -- (7, 4);
        \draw[blue, very thick] (0, 0) -- (2, 0) -- (2, 1) -- (4, 1) -- (4, 2) -- (6, 2) -- (6, 3) -- (7, 3) -- (7, 4);
    \end{scope}
    \end{tikzpicture}
    \end{center}
    \end{figure}
    
    \noindent Modified Christoffel word 
    \begin{equation*}
        \widetilde{\omega_{\frac{4}{7}}}(A, B) = ABABABB.
    \end{equation*}
    \noindent We have the IPM:
    
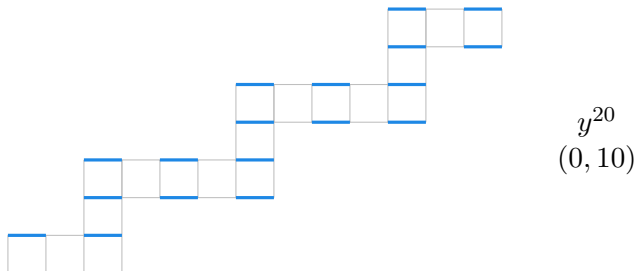
\begin{figure}[H]
    \begin{center}
    \begin{tikzpicture}[scale = 0.5]
        \draw[lightgray] (0, 0) grid (3, 1);
        \draw[lightgray] (2, 1) grid (3, 3);
        \draw[lightgray] (3, 2) grid (7, 3);
        \draw[lightgray] (6, 3) grid (7, 5);
        \draw[lightgray] (7, 4) grid (11, 5);
        \draw[lightgray] (10, 5) grid (11, 7);
        \draw[lightgray] (11, 6) grid (13, 7);

        \draw[blue, very thick] (0, 0) -- (1, 0);
        \draw[blue, very thick] (0, 1) -- (1, 1);
        \draw[blue, very thick] (2, 0) -- (3, 0);
        \draw[blue, very thick] (2, 1) -- (3, 1);
        \draw[blue, very thick] (2, 2) -- (3, 2);
        \draw[blue, very thick] (2, 3) -- (3, 3);
        \draw[blue, very thick] (4, 2) -- (5, 2);
        \draw[blue, very thick] (4, 3) -- (5, 3);
        \draw[blue, very thick] (6, 2) -- (7, 2);
        \draw[blue, very thick] (6, 3) -- (7, 3);
        \draw[blue, very thick] (6, 4) -- (7, 4);
        \draw[blue, very thick] (6, 5) -- (7, 5);
        \draw[blue, very thick] (8, 4) -- (9, 4);
        \draw[blue, very thick] (8, 5) -- (9, 5);
        \draw[blue, very thick] (10, 4) -- (11, 4);
        \draw[blue, very thick] (10, 5) -- (11, 5);
        \draw[blue, very thick] (10, 6) -- (11, 6);
        \draw[blue, very thick] (10, 7) -- (11, 7);
        \draw[blue, very thick] (12, 6) -- (13, 6);
        \draw[blue, very thick] (12, 7) -- (13, 7);

        \node at (15.5, 4) {$y^{20}$};
        \node at (15.5, 3) {$(0, 10)$};
    \end{tikzpicture}
    \end{center}
    \end{figure}
    \noindent Since $6 < b = 7$, we change the perfect matching in all vertical regions. The resulting perfect matching corresponds to the lattice point $(0, b) = (0, 7)$:
    
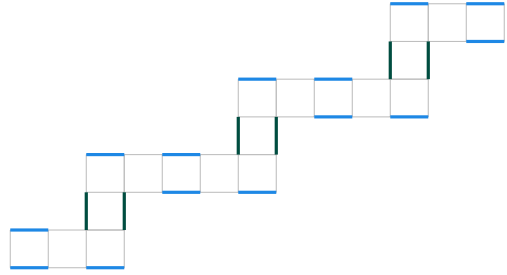
\begin{figure}[H]
    \begin{center}
    \begin{tikzpicture}[scale = 0.5]
        \draw[lightgray] (0, 0) grid (3, 1);
        \draw[lightgray] (2, 1) grid (3, 3);
        \draw[lightgray] (3, 2) grid (7, 3);
        \draw[lightgray] (6, 3) grid (7, 5);
        \draw[lightgray] (7, 4) grid (11, 5);
        \draw[lightgray] (10, 5) grid (11, 7);
        \draw[lightgray] (11, 6) grid (13, 7);

        \draw[blue, very thick] (0, 0) -- (1, 0);
        \draw[blue, very thick] (0, 1) -- (1, 1);
        \draw[blue, very thick] (2, 0) -- (3, 0);
        \draw[green, very thick] (2, 1) -- (2, 2);
        \draw[green, very thick] (3, 1) -- (3, 2);
        \draw[blue, very thick] (2, 3) -- (3, 3);
        \draw[blue, very thick] (4, 2) -- (5, 2);
        \draw[blue, very thick] (4, 3) -- (5, 3);
        \draw[blue, very thick] (6, 2) -- (7, 2);
        \draw[green, very thick] (6, 3) -- (6, 4);
        \draw[green, very thick] (7, 3) -- (7, 4);
        \draw[blue, very thick] (6, 5) -- (7, 5);
        \draw[blue, very thick] (8, 4) -- (9, 4);
        \draw[blue, very thick] (8, 5) -- (9, 5);
        \draw[blue, very thick] (10, 4) -- (11, 4);
        \draw[green, very thick] (10, 5) -- (10, 6);
        \draw[green, very thick] (11, 5) -- (11, 6);
        \draw[blue, very thick] (10, 7) -- (11, 7);
        \draw[blue, very thick] (12, 6) -- (13, 6);
        \draw[blue, very thick] (12, 7) -- (13, 7);

        \node at (15.5, 4) {$y^{14}z^6$};
        \node at (15.5, 3) {$(0, 7)$};
    \end{tikzpicture}
    \end{center}
    \end{figure}
    \noindent We are one diagonal above required, so we perform one operation on this intermediate perfect matching. Since there are not two boxes on the first line we cannot perform an addition and so perform a step $\Gamma_1^*$ (this can also be seen from the fact that the modified word begins with $ABA\dots$).
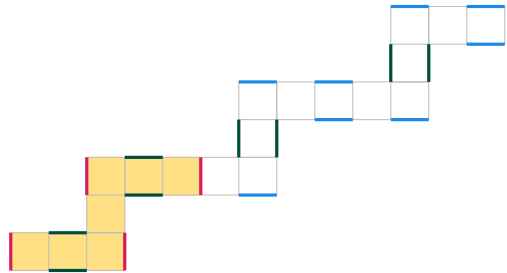
    \begin{figure}[H]
    \begin{center}
    \begin{tikzpicture}[scale = 0.5]
        \filldraw[yellow!50] (0, 0) -- (3, 0) -- (3, 2) -- (5, 2) -- (5, 3) -- (2, 3) -- (2, 1) -- (0, 1) -- (0, 0);
        \draw[lightgray] (0, 0) grid (3, 1);
        \draw[lightgray] (2, 1) grid (3, 3);
        \draw[lightgray] (3, 2) grid (7, 3);
        \draw[lightgray] (6, 3) grid (7, 5);
        \draw[lightgray] (7, 4) grid (11, 5);
        \draw[lightgray] (10, 5) grid (11, 7);
        \draw[lightgray] (11, 6) grid (13, 7);

        \draw[red, very thick] (0, 0) -- (0, 1);
        \draw[green, very thick] (1, 0) -- (2, 0);
        \draw[green, very thick] (1, 1) -- (2, 1);
        \draw[red, very thick] (3, 0) -- (3, 1);
        \draw[red, very thick] (2, 2) -- (2, 3);
        \draw[green, very thick] (3, 2) -- (4, 2);
        \draw[green, very thick] (3, 3) -- (4, 3);
        \draw[red, very thick] (5, 2) -- (5, 3);
        \draw[blue, very thick] (6, 2) -- (7, 2);
        \draw[green, very thick] (6, 3) -- (6, 4);
        \draw[green, very thick] (7, 3) -- (7, 4);
        \draw[blue, very thick] (6, 5) -- (7, 5);
        \draw[blue, very thick] (8, 4) -- (9, 4);
        \draw[blue, very thick] (8, 5) -- (9, 5);
        \draw[blue, very thick] (10, 4) -- (11, 4);
        \draw[green, very thick] (10, 5) -- (10, 6);
        \draw[green, very thick] (11, 5) -- (11, 6);
        \draw[blue, very thick] (10, 7) -- (11, 7);
        \draw[blue, very thick] (12, 6) -- (13, 6);
        \draw[blue, very thick] (12, 7) -- (13, 7);

        \node at (15.5, 4) {$x^4 y^8 z^8$};
        \node at (15.5, 3) {$(2, 4)$};
    \end{tikzpicture}
    \end{center}
    \end{figure}
    \noindent Now we are on the appropriate diagonal, but at the point $(2, 4)$. Therefore, we need to move two points along this diagonal in order to reach the required point.
    We continue through the snake, performing the required operations; firstly, a pull-back (since the next letter, in the modified word, after $ABA$ is a $B$) and then a twist (since the letter after that is an $A$).
    \begin{figure}[H]
    \begin{center}
    \begin{tikzpicture}[scale = 0.5]
        \filldraw[yellow!50] (0, 0) -- (3, 0) -- (3, 2) -- (7, 2) -- (7, 5) -- (6, 5) -- (6, 3) -- (2, 3) -- (2, 1) -- (0, 1) -- (0, 0);
        \draw[lightgray] (0, 0) grid (3, 1);
        \draw[lightgray] (2, 1) grid (3, 3);
        \draw[lightgray] (3, 2) grid (7, 3);
        \draw[lightgray] (6, 3) grid (7, 5);
        \draw[lightgray] (7, 4) grid (11, 5);
        \draw[lightgray] (10, 5) grid (11, 7);
        \draw[lightgray] (11, 6) grid (13, 7);

        \draw[red, very thick] (0, 0) -- (0, 1);
        \draw[green, very thick] (1, 0) -- (2, 0);
        \draw[green, very thick] (1, 1) -- (2, 1);
        \draw[red, very thick] (3, 0) -- (3, 1);
        \draw[red, very thick] (2, 2) -- (2, 3);
        \draw[green, very thick] (3, 2) -- (4, 2);
        \draw[green, very thick] (3, 3) -- (4, 3);
        \draw[green, very thick] (5, 2) -- (6, 2);
        \draw[green, very thick] (5, 3) -- (6, 3);
        \draw[red, very thick] (7, 2) -- (7, 3);
        \draw[red, very thick] (6, 4) -- (6, 5);
        \draw[red, very thick] (7, 4) -- (7, 5);
        \draw[blue, very thick] (8, 4) -- (9, 4);
        \draw[blue, very thick] (8, 5) -- (9, 5);
        \draw[blue, very thick] (10, 4) -- (11, 4);
        \draw[green, very thick] (10, 5) -- (10, 6);
        \draw[green, very thick] (11, 5) -- (11, 6);
        \draw[blue, very thick] (10, 7) -- (11, 7);
        \draw[blue, very thick] (12, 6) -- (13, 6);
        \draw[blue, very thick] (12, 7) -- (13, 7);

        \node at (15.5, 4) {$x^6 y^6 z^8$};
        \node at (15.5, 3) {$(3, 3)$};
    \end{tikzpicture}
    \end{center}
    \end{figure}
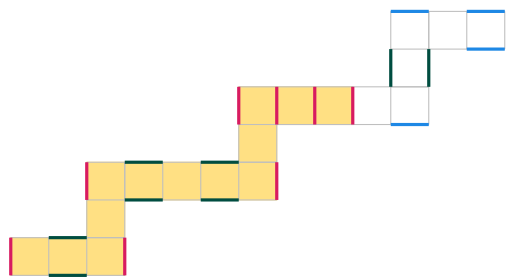
    \begin{figure}[H]
    \begin{center}
    \begin{tikzpicture}[scale = 0.5]
        \filldraw[yellow!50] (0, 0) -- (3, 0) -- (3, 2) -- (7, 2) -- (7, 4) -- (9, 4) -- (9, 5) -- (6, 5) -- (6, 3) -- (2, 3) -- (2, 1) -- (0, 1) -- (0, 0);
        \draw[lightgray] (0, 0) grid (3, 1);
        \draw[lightgray] (2, 1) grid (3, 3);
        \draw[lightgray] (3, 2) grid (7, 3);
        \draw[lightgray] (6, 3) grid (7, 5);
        \draw[lightgray] (7, 4) grid (11, 5);
        \draw[lightgray] (10, 5) grid (11, 7);
        \draw[lightgray] (11, 6) grid (13, 7);

        \draw[red, very thick] (0, 0) -- (0, 1);
        \draw[green, very thick] (1, 0) -- (2, 0);
        \draw[green, very thick] (1, 1) -- (2, 1);
        \draw[red, very thick] (3, 0) -- (3, 1);
        \draw[red, very thick] (2, 2) -- (2, 3);
        \draw[green, very thick] (3, 2) -- (4, 2);
        \draw[green, very thick] (3, 3) -- (4, 3);
        \draw[green, very thick] (5, 2) -- (6, 2);
        \draw[green, very thick] (5, 3) -- (6, 3);
        \draw[red, very thick] (7, 2) -- (7, 3);
        \draw[red, very thick] (6, 4) -- (6, 5);
        \draw[red, very thick] (7, 4) -- (7, 5);
        \draw[red, very thick] (8, 4) -- (8, 5);
        \draw[red, very thick] (9, 4) -- (9, 5);
        \draw[blue, very thick] (10, 4) -- (11, 4);
        \draw[green, very thick] (10, 5) -- (10, 6);
        \draw[green, very thick] (11, 5) -- (11, 6);
        \draw[blue, very thick] (10, 7) -- (11, 7);
        \draw[blue, very thick] (12, 6) -- (13, 6);
        \draw[blue, very thick] (12, 7) -- (13, 7);

        \node at (15.5, 4) {$x^8 y^4 z^8$};
        \node at (15.5, 3) {$(4, 2)$};
    \end{tikzpicture}
    \end{center}
    \end{figure}
    \noindent Now we have the perfect matching (at least one of them; there are $71$) corresponding to the lattice point $(4, 2)$.

    \begin{figure}[H]
    \begin{center}
    \begin{tikzpicture}
    \begin{scope}[xshift = 8cm, scale = 0.5]
        \filldraw[blue!25] (10, 0) -- (0, 10) -- (0, 7) -- (4, 0) -- (10, 0);
        \draw[lightgray] (0, 0) grid (10.5, 10.5);
        
        \draw[thick, ->] (0, 0) -- (0, 11);
        \draw[thick, ->] (0, 0) -- (11, 0);
        
        \node at (10, -0.5) {$10$};
        \node at (4, -0.5) {$4$};
        \node at (-0.7, 10) {$10$};
        \node at (-0.7, 7) {$7$};

        \node at (-0.5, 11) {$j$};
        \node at (11, -0.5) {$i$};

        \draw[orange] (0, 10) -- (0, 7);
        \draw[orange] (0, 7) -- (2, 4);
        \draw[purple] (2, 4) -- (4, 2);
        \filldraw[orange] (0, 10) circle (4pt);
        \filldraw[orange] (0, 9) circle (4pt);
        \filldraw[orange] (0, 8) circle (4pt);
        \filldraw[orange] (0, 7) circle (4pt);
        \filldraw[orange] (2, 4) circle (4pt);
        \filldraw[purple] (3, 3) circle (4pt);
        \filldraw[purple] (4, 2) circle (4pt);
    \end{scope}
    \end{tikzpicture}
    \caption{Path of perfect matchings from the above sequence of operations. Orange are left-most points, purple are from traversing the diagonal.}
    \label{fig:Example Path}
    \end{center}
    \end{figure}
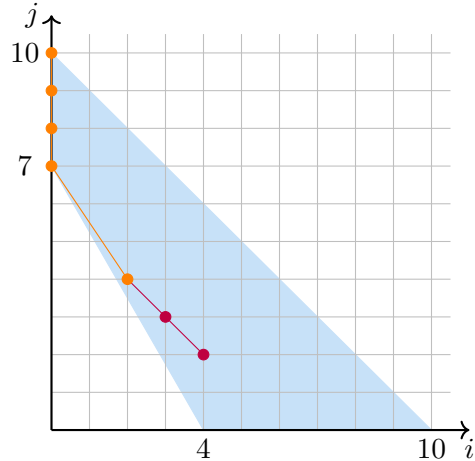
\end{example}

\section{Concluding remarks}
Now that we have seen that the combinatorial interpretation can be used to prove the saturation conjecture, it would be of interest to see if any further conclusions about the specific values of the coefficients can be drawn from taking a similar approach of walking along the Newton polygon with local operations on the perfect matchings of snake graphs.

\section*{Acknowledgements.}
I would like to thank Alexander Veselov and Brian Winn for helpful discussions about the formulation of the Saturation conjecture and initial ideas surrounding its proof. Additionally, Jason Semeraro for helpful discussions regarding the combinatorial approach. I also appreciate the support of the UK Research Council ESPRC for funding the project (EP/W028794/1).


\end{document}